  \documentclass{amsart}
 \usepackage{latexsym}
\usepackage{amssymb}
\usepackage{graphics}

\newtheorem{theorem}{Theorem}[section]
\newtheorem{corollary}{Corollary}[section]
\newtheorem{lemma}{Lemma}[section]

\newtheorem{remark}{Remark}[section]

\newtheorem{proposition}{Proposition}[section]

\def \a{\alpha }
\def \b {\beta}

\def \l{\lambda }

\begin{document}

\newcommand{\wta}{{\rm {wt} }  a }
\newcommand{\R}{\frak R}
\newcommand{\cV}{\mathcal V}
\newcommand{\cA}{\mathcal A}
\newcommand{\wtb}{{\rm {wt} }  b }
\newcommand{\bea}{\begin{eqnarray}}
\newcommand{\eea}{\end{eqnarray}}
\newcommand{\be}{\begin {equation}}
\newcommand{\ee}{\end{equation}}
\newcommand{\g}{\frak g}
\newcommand{\hg}{\hat {\frak g} }
\newcommand{\hn}{\hat {\frak n} }
\newcommand{\h}{\frak h}
\newcommand{\U}{\mathcal U}
\newcommand{\hh}{\hat {\frak h} }
\newcommand{\n}{\frak n}
\newcommand{\Z}{\Bbb Z}
\newcommand{\N}{{\Bbb Z} _{> 0} }
\newcommand{\Zp} {\Z _ {\ge 0} }
\newcommand{\C}{\Bbb C}
\newcommand{\Q}{\Bbb Q}
\newcommand{\1}{\bf 1}
\newcommand{\la}{\langle}
\newcommand{\ra}{\rangle}
\newcommand{\NS}{\bf{ns} }

\newcommand{\hf}{\mbox{$\tfrac{1}{2}$}}
\newcommand{\thf}{\mbox{$\tfrac{3}{2}$}}

\newcommand{\W}{\mathcal{W}}
\newcommand{\non}{\nonumber}
\def \l {\lambda}
\baselineskip=14pt
\newenvironment{demo}[1]%
{\vskip-\lastskip\medskip
  \noindent
  {\em #1.}\enspace
  }%
{\qed\par\medskip
  }

\def \l {\lambda}
\def \a {\alpha}

\keywords{vertex superalgebras, affine Lie algebras, Clifford
algebra, Weyl algebra, lattice vertex algebras,  critical level}
\title[]{
   Lie superalgebras and irreducibility of $A_1^{(1)}$--modules at the critical level }
\thanks{Partially supported by the MZOS grant 0037125 of the Republic of Croatia}
\author{ Dra\v zen Adamovi\' c }

\date{}
\curraddr{Department of Mathematics, University of Zagreb,
Bijeni\v cka 30, 10 000 Zagreb, Croatia} \email {adamovic@math.hr}
\markboth{Dra\v zen Adamovi\' c} { }
\bibliographystyle{amsalpha}
  \maketitle

\def \l {\lambda}
\def \a {\alpha}

\begin{abstract}
 We introduce the infinite-dimensional Lie superalgebra ${\cA}$
 and construct a family of  mappings  from certain category of ${\cA}$--modules to   the
 category of $A_1^{(1)}$--modules of critical level.  Using this
 approach, we prove the irreducibility of a  large family of
 $A_1^{(1)}$--modules at the critical level parameterized by ${\chi}(z) \in {\C}((z))$.
 As a consequence,  we present a new proof of irreducibility of certain Wakimoto
 modules.
 We also give a  natural realizations of irreducible quotients of
 relaxed Verma modules and calculate characters of these
 representations.
\end{abstract}

\section{Introduction}

In the analysis of certain Fock space representations of
infinite-dimensional Lie (super)algebras, one of the main problem
is to prove  irreducibility of these representations. Irreducible
highest weight representations of affine Lie algebras of critical
level can be realized  by using certain bosonic Fock
representation,
 called the Wakimoto modules (cf. \cite{W-mod}, \cite{FF1},
\cite{FB}, \cite{efren}, \cite{Scz}). Irreducibility of certain
Wakimoto modules gave a very natural proof of the Kac-Kazhdan
conjecture on characters of irreducible representations of
critical level (cf. \cite{KK}). On the other hand, the category of
representations of critical level is much richer than the category
$\mathcal{O}$. So one can investigate the modules outside the
category $\mathcal{O}$ and try to understand their structure.

 In particular one can investigate the relaxed Verma modules,
their irreducible quotients and the corresponding characters. Such
kind of representations appeared  in the context of representation
theory of the affine Lie algebra $\widehat{sl_2}$ on non-critical
levels (cf. \cite{FST}, \cite{AM}). In the present  paper we shall
demonstrate that these relaxed representations appear naturally at
the critical level and should be included in the representation
theory at this level. We will give a free field realization of the
irreducible quotients of relaxed Verma modules of critical level
and calculate their characters in the case of affine Lie algebra
$\widehat{sl_2}$.

Outside the critical level, the representation theory of affine
Lie algebra $A_1^{(1)}$ is related to the representation theory of
the $N=2$ superconformal algebra. In \cite{FST}, the authors
constructed mappings between certain categories of representations
of $\widehat{sl_2}$ and $N=2$ superconformal algebra. In the
context of vertex algebras these mappings was considered in
\cite{A3}. But at  the critical level the representation theory of
$\widehat{sl_2}$ is very different to the representation theory
outside the critical level. In particular, the associated vertex
algebra $N(-2\Lambda_0)$ contains  an infinite-dimensional center
(cf. \cite{efren}).

In the present paper we find an infinite-dimensional Lie
superalgebra ${\cA}$ with the important property that its
representation theory is related to those  of $\widehat{sl_2}$ at
the critical level. This algebra has  generators $G^{\pm} (r)$,
$T(n)$, $S(n)$, $r \in \hf + {\Z}$, $n \in {\Z}$, which satisfy
the following relations
\bea && [S(n), {\cA}] = [T(n), {\cA}]=0  \non \\&& \{ G ^{+} (r),
G ^{-} (s) \} = 2 S({r+s}) +
(r-s) T({r+s}) - ( r  ^ 2 - \frac{1}{4} ) \delta_{r+s,0} \non \\
 && \{ G ^{+} (r), G ^{+} (s) \}= \{ G ^{-} (r), G ^{-} (s) \} = 0 \non
\eea for all $n \in {\Z}$, $r,s \in {\hf} + {\Z}$.

The main difference between our algebra ${\cA}$ and the $N=2$
superconformal algebra is in the fact that ${\cA}$ contains large
center and that it doesn't contain the Virasoro and Heisenberg
subalgebra.  Next we consider the vertex superalgebra ${\cV}$
associated to a vacuum representation for ${\cA}$. This vertex
superalgebra is introduced in Section \ref{ver-def} as a vertex
subalgebra of $F \otimes M(0)$, where $F$ is a Clifford vertex
superalgebra and $M(0)$ a commutative vertex algebra.  Then
following \cite{A3} we show that there a non-trivial  vertex
algebra homomorphism $g: N(-2\Lambda_0) \rightarrow {\cV} \otimes
F_{-1}$, where $F_{-1}$ is a lattice vertex superalgebra
associated to the lattice ${\Z}{\b}$, $\la \b , \b \ra = -1$. This
result allows us to construct $N(-2\Lambda_0)$--modules from
${\cV}$--modules. Moreover, we prove that if $U$ is an irreducible
${\cV}$--module satisfying certain grading condition, then $U
\otimes F_{-1}=\oplus_{s\in {\Z} } {\mathcal L}_s (U)$ is a
completely reducible $N(-2\Lambda_0)$--module. Therefore every
component ${\mathcal L}_s (U)$ is an irreducible $A_1
^{(1)}$--module at the critical level.  It is important to notice
that the irreducibility result is proved by using the theory of
vertex algebras (cf. Lemma \ref{pom-ired}).

In this way the problem of constructing  irreducible $A_1
^{(1)}$--modules is reduced to the construction of irreducible
${\cV}$--modules. But on irreducible ${\cV}$--modules, the action
of the Lie superalgebra ${\cA}$ can be expressed by   the action
of generators of infinite-dimensional CLifford algebras.  By using
this fact, in Section \ref{konstr-a-ired} we prove the
irreducibility of a large family of ${\cV}$--modules. These
modules are parameterized by ${\chi}(z) \in {\C}((z))$. In Section
\ref{mappings} we construct mappings $\mathcal{L}_s$ which send
irreducible ${\cV}$--modules to the irreducible $\widehat{sl_2}$
modules at the critical level. As an application, in Section
\ref{wakimoto-mod} we present a proof of irreducibility of a large
family of the Wakimoto modules. In Section \ref{konstr-najv} we
study  the irreducible highest weight $\widehat{sl_2}$--modules of
critical level. In particular, we study the simple vertex algebra
$L(-2\Lambda_0)$. In Section \ref{konstr-relaxed} we get
realization of irreducible quotients of relaxed Verma modules. It
turns out that these irreducible modules can be realized on
certain  lattice type vertex algebra.

\section{ Vertex  algebra $N (k\Lambda_0)$ }

We make the assumption that the reader is familiar with the
axiomatic theory of vertex  superalgebras and  their
representations  (cf. \cite{DL}, \cite{FHL}, \cite{FLM},
\cite{LL}, \cite{K}, \cite{Z}).

In this section we recall some basic facts about vertex  algebras
associated to affine Lie algebras (cf. \cite{FZ}, \cite{Li},
\cite{MP}).

Let ${\g}$ be a finite-dimensional simple Lie algebra over ${\Bbb
C}$ and let $(\cdot,\cdot)$ be a nondegenerate symmetric bilinear
form on ${\g}$. Let  ${\g} = {\n}_- + {\h} + {\n}_+$    be a
triangular decomposition for ${\g}$.
% Let $\theta$ be the highest
%root for $\g$, end $e_{\theta}$ the corresponding root vector.
%Define $\rho$ as usual.
 The affine Lie algebra ${\hg}$ associated
with ${\g}$ is defined as $ {\g} \otimes {\C}[t,t^{-1}] \oplus
{\C}c   $ where $c$ is the canonical central element \cite{K-b}
 and  the Lie algebra structure
is given by $$ [ x \otimes t^n, y \otimes t^m] = [x,y] \otimes t
^{n+m} + n (x,y) \delta_{n+m,0} c.$$    We will write $x(n)$ for
$x \otimes t^{n}$.

%Adjoining the degree operator $d$ such that $[d, x(n)]=n x(n)$ and
%$[d,c]=0$ gives us the Lie algebra $ \widetilde{\g}$, the affine
%Kac-Moody Lie algebra $A_1^{(1)}$.

The Cartan subalgebra ${\hh}$ and  subalgebras ${\hg}_+$,
${\hg}_-$ of ${\hg}$ are defined by $${\hh} = {\h} \oplus {\C}c ,
\quad
 {\hg}_{\pm} =  {\g}\otimes t^{\pm1} {\C}[t^{\pm1}].$$

Let  $P = {\g}\otimes {\C}[t] \oplus {\C}c $ be upper parabolic
subalgebra.   For every $k \in {\C}$,  let ${\C} v_k$ be
$1$--dimensional $P$--module  such that the subalgebra ${\g}
\otimes  {\C}[t] $ acts trivially,
 and  the central element
$c$ acts as multiplication with $k \in {\C}$. Define the
generalized Verma module $N( k \Lambda_0)$ as
$$N({k}\Lambda_0) = U(\hg) \otimes _{ U(P) } {\C} v_k .$$
Then $ N(k \Lambda_0)$ has a natural structure of a vertex
algebra. The vacuum vector is ${\1} = 1 \otimes v_k$.

The vertex algebra $N(k\Lambda_0)$ has very rich representation
theory. Let $U$ be any ${\g}$--module. Then  $U$ can be consider
as a $P$--module. The induced ${\hg}$--module $N(k, U) = U({\hg})
\otimes_{ U(P)} U $ is a module for the vertex algebra
$N(k\Lambda_0)$.

Let $N^{1}(k \Lambda_0)$ be the maximal ideal in the vertex
algebra $N(k \Lambda_0)$. Then $L(k \Lambda_0) = \frac{N(k
\Lambda_0)}{N^{1}( k \Lambda_0)}$ is a simple vertex  algebra.

 Let now $\g =sl_2(\C) $ with generators $e$, $f$, $h$ and
relations $[h,e]= 2 e$, $[h, f] = -2 f$, $[e,f]= h$. Let
$\Lambda_0$, $\Lambda_1$  be the fundamental weights for $\hg$.

For $s \in {\Z}$, we define $H^{s} = - s \frac{h}{2}$. Then $(
H^{s} , h ) = -s$.

Define $$\Delta_s (z) = z^{ H^{s}(0)} \exp \left(\sum_{ n=1}
^{\infty} \frac{H^{s}(n)}{-n} (-z)^{-n}\right).$$
Applying  the results obtained in  \cite{Li5} on
$N(k\Lambda_0)$--modules we get the following proposition.

\begin{proposition} \label{novi-op} Let $s \in {\Z}$.
%
% \item[(1)]
For any  $N(k\Lambda_0)$--module $(M,Y_M(\cdot,
 z))$,
$$ (\pi_s(M), Y_M ^{s}(\cdot,z)):=(M, Y_M (\Delta(H^{s}, z) \cdot, z))$$
is a  $N(k \Lambda_0)$--module. $\pi_s(M)$ is an irreducible weak
$N(k\Lambda_0)$--module if and only if
  $M$ is   an irreducible weak $N(k\Lambda_0)$--module.
\end{proposition}

By definition we have:
\bea &&  \Delta(H^{s}, z) e(-1){\1} = z^{-s} e(-1){\1}, \nonumber \\
 &&
\Delta(H^{s} ,z) f(-1){\1} = z^{s} f(-1){\1}, \nonumber  \\
&& \Delta(H^{s}, z) h(-1){\1} = h(-1){\1} - s k z^{-1} {\1}.
\nonumber \eea

In other words, the corresponding automorphism $\pi_s$ of $U(\hg)$
satisfies the condition:
$$ \pi_s ( e(n) ) = e(n - s), \ \ \pi_s ( f (n) ) = f (n+s), \ \
\pi_s( h(n) ) = h(n) - s k \delta_{n,0}.$$

In the case $s=-1$, one can see  that  $$\pi_{-1}( L((k-n)
\Lambda_0 + n\Lambda_1)) = L(n \Lambda_0 + (k-n) \Lambda_1), \quad
\mbox{for every} \ n \in {\Zp}.$$
It is also important to notice the following important property:
$$ \pi_{s + t} (M) \cong \pi_s ( \pi_t (M) ), \ \ (s, t \in
{\Z}).$$
In particular, $$ M \cong \pi_0 (M) \cong \pi _s ( \pi_{-s} (M) ).
$$

\begin{lemma} \label{sing-za}
 Let $x \in {\C}$.  Assume that $U$ is an irreducible
$N(k\Lambda_0)$--module which is generated by the vector $v_{s}$
$(s \in {\Z})$ such that:
\bea && e(n- s) v_{s} = f (n+s+1) v_{s}=0 \quad (n \ge 0),
\label{rel-0-1} \\
&& h(n) v_{s} = \delta_{n,0} ( x +  k s) v_{s} \quad (n \ge
0).\label{rel-0-2} \eea
Then
$$ U \cong \pi_{-s} ( L((k-x)
\Lambda_0 + x \Lambda_1). $$
\end{lemma}
{\em Proof.}  We consider the   $N(k\Lambda_0)$--module $\pi_{s}
(U)$. By construction, we have that $\pi_{s}(U)$ is an irreducible
highest weight $\hg$--module with the highest weight $(k-x)
\Lambda_0 + x \Lambda_1 $. Therefore,
$\pi_{s}(U) \cong L( (k-x) \Lambda_0 + x \Lambda_1),$
which implies that
$$ U \cong \pi_{-s} ( \pi_{s} (U ) ) \cong \pi_{-s} ( L( (k-x) \Lambda_0 + x \Lambda_1)  ),$$
and the Lemma  holds.    \qed

\section{ Clifford vertex superalgebras}

   The Clifford
algebra $CL$ is a complex associative algebra generated by
$$ \Psi^{\pm}(r) , \  r  \in \hf + {\Z},$$ and relations
\bea
&& \{\Psi^{\pm}(r) , \Psi^{\mp}(s) \} = \delta_{r+s,0}; \quad
 \{\Psi^{\pm}(r) , \Psi^{\pm}(s)\}=0
\nonumber
\eea
where $r, s \in {\hf}+ {\Z}$.

 Let $F$ be the irreducible $CL$--module generated by
 the
cyclic vector ${\1}$ such that
$$ \Psi^{\pm} (r) {\1} = 0 \quad
\mbox{for} \ \ r > 0 .$$
A basis of $F$ is given by
$$ \Psi^{+ }({-n_1-{\hf}})  \cdots \Psi^{+}({-n_r-{\hf}})  \Psi^{-}({-k_1-{\hf}})  \cdots \Psi^{-}({-k_s-{\hf}})
 {\1} $$
where $n_i, k_i \in {\Zp}$,  $n_1 >n_2 >\cdots >n_r  $, $k_1 >k_2
>\cdots
>k_s $.

Define the following   fields on $F$
$$   \Psi^{+}(z) = \sum_{ n \in   {\Z}
 } \Psi^{+}(n+{\hf} )  z ^{-n- 1}, \quad  \Psi^{-} (z) = \sum_{ n \in {\Z}
 } \Psi^{-} (n+{\hf} )  z ^{-n-1}.$$

 The fields $\Psi^{+}(z)$ and $\Psi^{-}(z)$ generate on $F$  the
unique structure of a simple vertex superalgebra (cf. \cite{Li},
\cite{K}, \cite{FB}).

Define the following Virasoro vector in $F$ :

$$ \omega ^{(f)} = \frac{1}{2} ( \Psi ^{+} (-\tfrac{3}{2} ) \Psi
^{-}(-\tfrac{1}{2}) +  \Psi ^{-} (-\tfrac{3}{2} ) \Psi
^{+}(-\tfrac{1}{2})) {\1}.$$

Then the components of the field  $L^{(f)}(z) = Y(\omega^{(f)},z)
= \sum_{n \in {\Z} } L^{(f)}(n) z^{-n-2}$ defines on $F$ a
representation of the Virasoro algebra with central charge
$c^{(f)}=1$.

Set $$J^{(f)}(z)= Y(\Psi^{+}(-{\hf})\Psi^{-}(-{\hf}){\1},z)=
\sum_{n \in {\Z} } J^{(f)}(n) z^{-n-1}.$$

Then we have
$$ [J^{(f)}(n), \Psi^{\pm} (m+ {\hf})] = \pm \Psi^{\pm} (m+n +
{\hf}).$$

 Let $\widetilde{F} = \mbox{Ker}_F \Psi ^{-}(\hf)$
be the subalgebra of the vertex superalgebra $F$ generated by the
fields
$$ \partial \Psi^{+}(z)=\sum_{n \in {\Z} } -n \Psi^{+}(n - {\hf}) z^{-n-1} \ \ \mbox{and} \ \
 \Psi^{-}(z) = \sum_{ n \in {\Z} }   \Psi^{-}(n+ {\hf})
z^{-n-1}. $$
Then $\widetilde{F}$ is a simple vertex superalgebra with basis
\bea \label{baza-tilde} && \Psi^{+ }({-n_1-{\hf}})  \cdots
\Psi^{+}({-n_r-{\hf}}) \Psi^{-}({-k_1-{\hf}})  \cdots
\Psi^{-}({-k_s-{\hf}})
 {\1} \eea
where $n_i, k_i \in {\Zp}$,  $n_1 >n_2 >\cdots >n_r \ge 1 $, $k_1
>k_2
>\cdots
>k_s \ge 0 $.

 Let $\overline{F}= \mbox{Ker}_F \Psi ^{-}(\hf) \cap \mbox{Ker}_F \Psi ^{+}(\hf)$
be the subalgebra of the vertex superalgebra $F$ generated by the
fields
$$ \partial \Psi^{+}(z)=\sum_{n \in {\Z} } -n \Psi^{+}(n - {\hf}) z^{-n-1} \ \ \mbox{and} \ \
\partial \Psi^{-}(z) = \sum_{ n \in {\Z} } - n  \Psi^{-}(n-{\hf}) z^{-n-1}. $$
Then $\overline{F}$ is a simple vertex superalgebra with basis
\bea \label{baza-bar} && \Psi^{+ }({-n_1-{\thf}})  \cdots
\Psi^{+}({-n_r-{\thf}}) \Psi^{-}({-k_1-{\thf}})  \cdots
\Psi^{-}({-k_s-{\thf}})
 {\1} \eea
where $n_i, k_i \in {\Zp}$,  $n_1 >n_2 >\cdots >n_r  $, $k_1 >k_2
>\cdots
>k_s $.

\section{ The vertex superalgebra ${\cV}$}
\label{ver-def}
In this section we shall define the vertex superalgebra ${\cV}$
and study its representation theory. The vertex superalgebra
${\cV}$  contains a large center. Moreover, the  vertex
superalgebra $\overline{F}$ is a simple quotient of ${\cV}$.

Let $M(0) = {\C}[{\gamma}^{+}(-n), {\gamma}^{-}(n) \  \vert \ n <0
]$ be the commutative vertex algebra generated by the fields
$${\gamma} ^{\pm} (z) = \sum_{ n  < 0} {\gamma}^{\pm} (n) ^{-n-1}. $$
(cf. \cite{efren}). Let $\chi^{\pm}(z) = \sum_{n \in {\Z} }
\chi^{\pm}_{n} z^{-n-1} \in {\C} ((z))$. Let $M(0, \chi^{+},
\chi^{-})$ denotes the $1$--dimensional irreducible $M(0)$--module
with the property that every element ${\gamma}^{\pm}(n)$ acts on
$M(0, \chi^{+}, \chi^{-})$ as multiplication with $\chi^{\pm}_n
\in {\C}$.

Let now ${\mathcal F}$ be the vertex superalgebra generated by the
fields $\Psi^{\pm} (z)$ and ${\gamma}^{\pm}(z)$. Therefore
${\mathcal F} = F \otimes M(0)$. Denote by  ${\cV}$ the vertex
subalgebra of the  vertex superalgebra ${\mathcal F}$  generated
by the following vectors
\bea
 \tau^{\pm} &=& (\Psi^{\pm} (-\tfrac{3}{2}) + {\gamma}^{\pm} (-1)
 \Psi^{\pm}(-\hf)) {\1}, \label{def-tau} \\
 j &=& \frac{ {\gamma}^{+} (-1) - {\gamma}^{-}(-1)}{2} {\1}, \label{def-j} \\
 \nu &=& \frac{ 2 {\gamma}^{+} (-1) {\gamma}^{-}(-1) + {\gamma}^{+}(-2) + {\gamma}^{-}(-2)}{4}
  {\1} .
\label{def-nu}  \eea
Then the  vertex superalgebra structure on ${\cV}$     is
generated by the following fields \bea && G^{\pm} (z) = Y(\tau
^{\pm} ,z)
= \sum _{n \in {\Z} } G ^{\pm} (n+ {\hf} ) z ^{-n-2}, \non \\
 && S (z) = Y(\nu,z)
= \sum _{n \in {\Z} } S({n })  z ^{-n-2}, \non \\
&& T (z) = Y(j,z) = \sum _{n \in {\Z} } T(n )  z ^{-n-1}. \non
\eea
By using commutator formulae, we have that  the components of
these fields span an infinite-dimensional Lie superalgebra. Let us
denote  this Lie superalgebra by ${\cA}$. This algebra has
generators $G^{\pm} (r)$, $T(n)$, $S(n)$, $r \in \hf + {\Z}$, $n
\in {\Z}$, which satisfy the following relations
\bea && [S(n), {\cA}] = [T(n), {\cA}]=0 ,  \non \\&& \{ G ^{+}
(r), G ^{-} (s) \} = 2 S({r+s}) +
(r-s) T({r+s}) - ( r  ^ 2 - \frac{1}{4} ) \delta_{r+s,0} , \non \\
 && \{ G ^{+} (r), G ^{+} (s) \}= \{ G ^{-} (r), G ^{-} (s) \} = 0 \non
\eea for all $n \in {\Z}$, $r,s \in {\hf} + {\Z}$.

So the vertex superalgebra ${\cV}$ is generated by the Lie
superalgebra ${\cA}$. Thus we can study ${\cV}$--modules as
modules for the Lie superalgebra ${\cA}$. The proof of the
following proposition is standard.
\begin{proposition} We have:
\item[(1)] ${\cV} = U(\cA) . {\1}$
\item[(2)] Assume that $U$ is a ${\cV}$--module. Then $U$ is an
irreducible ${\cV}$--module if and only if $U$ is an irreducible
${\cA}$--module.
\end{proposition}

Let ${\cV} ^{com}$ be the vertex subalgebra of ${\cV}$ generated
by the fields $S(z)$ and $T(z)$. ${\cV} ^{com}$ is a commutative
vertex algebra.

The operator $J^{f}(0)$ acts semisimply on the vertex superalgebra
${\cV}$ and defines the following ${\Z}$--graduation:
\bea  {\cV} &=& \bigoplus _{m \in {\Z} } {\cV} ^{m}, \quad
\mbox{where}
\non \\
{\cV} ^{m} &=& \{ v \in {\cV} \ \vert \ J^{f}(0) v = m v \}
\non \\
&=& \mbox{span}_{\C} \{ G^{+ }({-n_1-{\thf}})  \cdots
G^{+}({-n_r-{\thf}})  G^{-}({-k_1-{\thf}})  \cdots
G^{-}({-k_s-{\thf}}) w  \vert \non \\
&& w \in {\cV} ^{com}, n_i, k_j \in {\Zp}, r-s=m \}.
\label{opis-komponente} \eea

It is clear that ${\cV} ^{com} \subset {\cV} ^{0}$.

Every $U({\cA})$-submodule of ${\cV}$ becomes an ideal in the
vertex superalgebra ${\cV}$. Let $I^{com}= U({\cA}). {\cV} ^{com}$
be the ideal in ${\cV}$ generated by ${\cV} ^{com}$. From the
definition of vertex superalgebras ${\cV}$ and $\overline{F}$ we
get the following result.

\begin{proposition} \label{simple-quotient} The quotient vertex superalgebra ${\cV} / I^{com}  $ is
isomorphic to the simple vertex superalgebra $\overline{F}$.
\end{proposition}

\section{Irreducibility of certain ${\cV}$--modules}
 \label{konstr-a-ired}

 In this section  we shall consider a family of
 irreducible ${\cV}$--modules.

For ${\chi}^{+}, {\chi}^{-} \in {\C}((z))$ we set
$F({\chi}^{+}, {\chi}^{-}) :=F \otimes M(0,\chi^{+},\chi^{-})$.

Then $F({\chi}^{+}, {\chi}^{-})$ is a module for the vertex
superalgebra ${\cV}$, and therefore for the Lie superalgebra
${\cA}$.

Since $ M(0,\chi^{+},\chi^{-})$ is one-dimensional, we have that
as  a vector space

\bea  \label{identification-prva} F({\chi}^{+}, {\chi}^{-}) \cong
F \cong { \bigwedge} (\Psi^{\pm} (-i-\hf) \ \vert \ i \ge 0 ).
\eea

This actually shows that for every ${\chi} ^{+}, {\chi}^{-} \in
{\C}((z))$ on the vertex superalgebra $F$ exists  the structure of
a ${\cA}$--structure. In this section we shall use this
identification.

\begin{proposition} \label{vazna-1}
  Assume that ${\l} \in {\C} \setminus {\Z}$, $p \in {\Zp}$ and
that
$$ {\chi}(z) = \sum_{n=-p} ^{\infty} {\chi}_{-n} z ^{n-1} \in {\C}((z))$$
satisfies the following conditions
\bea
&& {\chi}_p \ne 0, \ \label{uvjet-pr} \\
&& {\chi}_0 \in   {\C} \setminus {\Z}  \ \ \mbox{if} \ p=0
\label{uvjet-dr}. \eea
Then $F(\frac{\l}{z}, {\chi} )$ is an irreducible ${\cV}$--module.
\end{proposition}
{\em Proof.} Since $F(\frac{\l}{z},{\chi}) $ is a ${\cV}$--module,
it remains to prove that it is an irreducible module for the Lie
superalgebra ${\cA}$. Recall (\ref{identification-prva}).  The
${\cA}$--module structure on $F(\frac{\l}{z}, {\chi})$ is uniquely
determined by the following action of generators of ${\cA}$ on
$F$:
\bea
G^{+}(i-\hf) &=& ({\l} -i)  \Psi^{+}(i-\hf) \label{+djelovanje} , \\
G^{-}(i-\hf) &=&- i \Psi^{-}(i-\hf) + \sum_{k=-p} ^{\infty}
{\chi}_{-k} \Psi^{-} (k+i -{\hf}). \label{-djelovanje}
 \eea
 First we shall prove that the vacuum vector is a cyclic vector of
 the $U(\cA)$--action,i.e.,
 \bea
\label{ciklic} U({\cA}). {\1} = F.
 \eea

Take an arbitrary  basis element
\bea &&v= \Psi^{+ }({-n_1-{\hf}})  \cdots \Psi^{+}({-n_r-{\hf}})
\Psi^{-}({-k_1-{\hf}})  \cdots \Psi^{-}({-k_s-{\hf}})
 {\1}  \in F ,  \label{baza}\eea
 where $n_i, k_i \in {\Zp}$,  $n_1 >n_2 >\cdots
>n_r \ge 0 $, $k_1
>k_2
>\cdots
>k_s \ge 0 $.

Let $N \in {\Zp}$ such that $N \ge k_1$. By using
(\ref{-djelovanje}) we get that
$$ G^{-}(p- N - {\hf}) \cdots G^{-}(p-\thf) G^{-}(p-\hf) {\1} = C \Psi^{-}(-N -{\hf}) \cdots
\Psi ^{-}(-\thf)\Psi^{-} (-\hf) {\1},
$$
where
$$C= \left\{ \begin{array}{cc}
   {\chi}_p ^{N+1} & \mbox{if} \ p \ge 1  \\
  {\chi}_0 ({\chi}_0 +1) \cdots ({\chi}_0 + N) & \ \mbox{if} \ p=0
\end{array} . \right. $$
So $C\ne 0$, and we have that
$$
\Psi^{-}(-N -{\hf}) \cdots \Psi ^{-}(-\thf)\Psi^{-} (-\hf) {\1}
\in U({\cA}) .{\1} .$$
By using this fact and the action of elements $G^{+}(i-{\hf})$, $i
\in {\Z}$, we obtain that $v \in U(\cA).{\1}$. In this way  we
proved (\ref{ciklic}).

In order to prove irreducibility, it is enough to show that
arbitrary basis  element $v$ of the form (\ref{baza}) is cyclic in
$F$.  This follows from (\ref{ciklic}) and the fact that
$$
G^{-}(n_r+p+\hf) \cdots G^{-}(n_1+p+\hf) G^{+}(k_s+\hf) \cdots
G^{+}(k_1+\hf) v = C' {\1},
$$
where the non-trivial constant $C'$ is given by
\bea && C' = (-1) ^{r + s} ({\l}-k_1 -1) \cdots ({\l} - k_s -1)
\cdot C'' \ \ \mbox{and} \ \nonumber \\   && C''= \left\{
\begin{array}{cc}
   {\chi}_p ^{r} & \mbox{if} \ p \ge 1  \\
  ({\chi}_0 -n_1-1) \cdots ({\chi}_0 -n_r-1) & \ \mbox{if} \ p=0
\end{array} . \right. \nonumber \eea
\qed

This proposition has the following important consequence.
\begin{corollary} \label{ireducibilnost-prva}
 Assume that $\l, \mu \in {\C} \setminus
{\Z}$. Then $F(\frac{\l}{z},\frac{\mu}{z})$ is an irreducible
${\cV}$--module.
\end{corollary}

\vskip 5mm

Now let ${\chi} (z)  \in {\C}((z))$. Define :
$$ \widetilde{F}_{\chi} := \widetilde{F} \otimes M(0,0,{\chi}).$$

It is clear that $\widetilde{F}_{\chi}$ is a submodule of the
${\cV}$--module $F ( 0,{\chi})$. Now we shall prove the following
important irreducibility result:

\begin{proposition} \label{vazna-2}
Assume that $p \in {\Zp}$ and that
$$ {\chi}(z) = \sum_{n=-p} ^{\infty} {\chi}_{-n} z ^{n-1} \in {\C}((z))$$
satisfies the following conditions
\bea
&& {\chi}_p \ne 0,   \label{uvjet-prvi} \\
&& {\chi}_0 \in  \{1\} \cup \left({\C} \setminus {\Z}\right)  \ \
\mbox{if} \ p=0 \label{uvjet-drugi}. \eea
Then
$\widetilde{F}_{\chi}$ is an irreducible ${\cV}$--module.
\end{proposition}
{\em Proof.} Since $\widetilde{F}_{\chi}$ is a ${\cV}$--module, it
remains to prove that $\widetilde{F}_{\chi}$ is an irreducible
module for the Lie superalgebra ${\cA}$. The ${\cA}$--module
structure on $\widetilde{F}_{\chi}$ is uniquely determined by the
following action of the Lie superalgebra ${\cA}$ on
$\widetilde{F}$:
\bea
G^{+}(i-\hf) &=& -i \Psi^{+}(i-\hf) \label{++djelovanje} , \\
G^{-}(i-\hf) &=&- i \Psi^{-}(i-\hf) + \sum_{k=-p} ^{\infty}
{\chi}_{-k} \Psi^{-} (k+i -{\hf}). \label{--djelovanje}
 \eea
 By using this action, the basis description   (\ref{baza-tilde}) of $\widetilde{F}$, and   same proof to those  of
 Proposition \ref{vazna-1} we get  the irreducibility result.
 \qed

\begin{corollary} \label{ireducibilnost-treca}
Assume that $\l \in {\C} \setminus {\Z}$ or $\l = -1$. Then
$ \widetilde{F}_{-\tfrac{\l}{z} }$ is an irreducible
${\cV}$--module.
\end{corollary}

\begin{proposition} \label{ireducibilnost-druga}
 Assume that $m, n \in {\Zp}$. Then
\bea
\overline{F}(-\tfrac{m}{z},-\tfrac{n}{z}) &=& \mbox{Ker}_{
F(-\tfrac{m}{z},-\tfrac{n}{z})} G^{-}(m+\hf) \bigcap
\mbox{Ker}_{ F(-\tfrac{m}{z},-\tfrac{n}{z})} G^{+}(n+\hf) \nonumber \\
&\cong& \bigwedge (\Psi^{+}(-i-\hf),  \Psi^{-}(-j-\hf) \ \vert \
i, j \in {\Zp}, \ i\ne m, \ j \ne n ) \nonumber
\eea
is an irreducible ${\cV}$--module. In particular, the
${\cV}$--module $\overline{F}(0,-\tfrac{n}{z})$ is irreducible.
\end{proposition}
{\em Proof.} First we notice that on
$F(-\frac{m}{z},-\frac{n}{z})$
\bea \label{non-str} && G^{+} (i-{\hf}) =- (i+m) \Psi^{+} (i -\hf)
, \quad G^{-} (i-{\hf}) = - (i+n) \Psi^{-} (i -\hf) . \eea
This implies that
\bea \label{cl-mn} &&\overline{F}(-\frac{m}{z},-\frac{n}{z}) \cong
\bigwedge (\Psi^{+}(-i-\hf), \Psi^{-}(-j-\hf) \ \vert \ i, j \in
{\Zp}, \ i\ne m, \ j \ne n )\eea
is a ${\cV}$--submodule of $F(-\frac{m}{z},-\frac{n}{z})$. It
remains to prove that $\overline{F}(-\frac{m}{z},-\frac{n}{z})$ is
an irreducible ${\cA}$--module. By using (\ref{non-str}) we have
that   ${\overline F}(-\frac{m}{z},-\frac{n}{z})$ has the
structure of a module for the subalgebra $ ^{m,n} CL$ of the
Clifford algebra $CL$ generated by
$$ \Psi^{+} (-i-\hf), \Psi^{-}(-j-\hf), \ i , j \in {\Z}, \ i \ne
m, \ j \ne n. $$
Since $\bigwedge (\Psi^{+}(-i-\hf), \Psi^{-}(-j-\hf) \ \vert \ i,
j \in {\Zp}, \ i\ne m, \ j \ne n )$ is an irreducible $^{m,n}
CL$--module, we have that $\overline{F}
(-\frac{m}{z},-\frac{n}{z})$ is an irreducible ${\cA}$--module.
\qed

\section{ Vertex algebras at the critical level }
\label{mappings}

 In previous sections we investigated the
properties of the vertex superalgebra ${\cV}$. This vertex
superalgebra is generated by the Lie superalgebra ${\cA}$ which is
similar to the $N=2$ superconformal algebras. This makes the
vertex superalgebra ${\cV}$ similar to the  $N=2$ vertex
superalgebras investigated in \cite{EG}, \cite{FST} and \cite{A3}.
 The main difference is that ${\cV}$ doesn't contain the
Virasoro and the Heisenberg subalgebra. On the other hand ${\cV}$
 contains large center.
  One important property of the $N=2$
vertex superalgebras is their connection to the affine
$\widehat{sl_2}$--vertex algebras. The most effective way for
studying this connection is by using Kazama-Suzuki and anti
Kazama-Suzuki mapping(cf. \cite{KS}, \cite{FST}). Motivated by the
anti Kazama-Suzuki mapping, we shall get a realization of the
vertex algebras  associated to the   $\widehat{sl_2}$--modules at
the critical level.

Let $F_{-1} $,  be the lattice vertex superalgebra $V_{ \Z \b}$,
associated to the lattice ${\Z} \b$, where $\la \b, \b \ra = -1$.
As a vector space, $F_{-1}$ is isomorphic to $M_{\b}(1)\otimes
{\C}[L]$, where $M_{\b}(1)$ is a level one irreducible module for
the Heisenberg algebra $\hat{\h}_{\Z}$ associated to the one
dimensional abelian algebra ${\h}= L \otimes_{\Z} {\C}$ and
${\C}[L]$ is the group algebra with a generator $e^{\b}$.  The
generators of $F_{-1}$ are $e^{\b}$ and $e^{-\b}$. Moreover,
$F_{-1}$ is a simple vertex superalgebra and a completely
reducible $M_{\b} (1)$--module isomorphic to
$$F_{-1}\cong \bigoplus_{m \in {\Z} } F_{-1} ^{m},$$ where $ F_{-1} ^{m}$ is
an irreducible $M_{\b}(1)$--module generated by $e^{m\b}$.

 We
shall now consider the vertex superalgebra ${\cV} \otimes F_{-1}$.
Let $Y$ be the  vertex operator defining the vertex operator
superalgebra structure
  on ${\cV} \otimes F_{-1}$.
  For every $v \in {\cV} \otimes
  F_{-1}$, let $Y(v,z) = \sum_{ n \in {\Z} } v_n z^{-n-1}$.

 Define
\bea  && e= G^{+}(-\tfrac{3}{2}) {\1} \otimes e^{\b} ,
\label{defin-e} \\  && h= - 2 ( {\1} \otimes \b(-1)  - T(-1)
\otimes {\1}), \label{defin-h} \\ && f= G^{-} (-\tfrac{3}{2}) {\1}
\otimes e^{-\b}. \label{defin-f} \eea
For $x \in \mbox{span}_{\C} \{e,f,h\}$ set $x(z) = Y(x,z)= \sum_{n
\in {\Z} } x(n) z^{-n-1}$. Then the components of the field
$e(z)$, $f(z)$ and $h(z)$ satisfy the commutation relations for
the affine Lie algebra $\hat{sl_2}$ of level $-2$.
In particular we have
\bea
&&e(n) = \sum_{i \in {\Z} }   G^{+}(i - {\hf} ) \otimes
e^{\b}_{n-i-1}, \nonumber
\\
&& h(n)= - 2  \b(n) + 2 T(n) , \nonumber \\
&& f(n) =  \sum_{i \in {\Z} }   G^{-}(i-{\hf} ) \otimes
e^{-\b}_{n-i-1}. \nonumber \eea

 Denote by $V$
the subalgebra  of ${\cV} \otimes F_{-1}$ generated by $e$, $f$
and $h$, i.e.,
$$  V=\mbox{span}_{\C} \{ u^{1}_{n_1} \cdots u^{r} _{n_r}
({\1} \otimes {\1})  \vert \   u^{1}, \dots, u^{r} \in \{e,f,h\},
\ n_1, \dots, n_r \in {\Z}, r \in {\Zp} \}. $$
As a $U(\hg)$--module, $V$ is a cyclic module generated by the
vacuum vector ${\1}\otimes {\1}$. This implies that $V$ is a
certain quotient of the vertex algebra $N(-2 \Lambda_0)$.

So we have:

\begin{proposition}
There exists a non-trivial homomorphism of   vertex algebras  $$
g: N(-2\Lambda_0) \rightarrow {\cV} \otimes F_{-1} , $$
which is uniquely determined by (\ref{defin-e})-(\ref{defin-f}).
\end{proposition}

 \vskip 3mm

Define the vector
\bea \label{def-operatora-h} && H = \Psi^{+}(-\hf) \Psi^{-} (-\hf)
{\1} \otimes {\1} + {\1} \otimes \b (-1) \in F \otimes F_{-1}
\subset {\mathcal F} \otimes F_{-1}. \eea
 Then the operator $H(0) $
acts semisimply on vertex superalgebras ${\mathcal F} \otimes
F_{-1}$ and ${\cV} \otimes F_{-1}$. Let
$$
   W(s) = \{ v \in {\cV} \otimes F_{-1} \  \vert H(0) v = s v
\} .
$$
Then we have the following decomposition:
$$
{\cV} \otimes F_{-1} \cong \oplus_{s \in {\Z} } W(s).$$

\vskip 5mm

  Let ${\hat t}$
be the (commutative) Lie algebra generated by the components of
the field $T(z)$, and  let $M_T (0) \cong {\C}[T(-1), T(-2),
\dots]$ be the associated (commutative) vertex algebra.

Let ${\hg} ^{ext} = {\hg} \oplus {\hat t}$ be the extension of the
affine Lie algebra ${\hg}$ by the Lie algebra ${\hat t}$ such that
every element of ${\hat t}$ is in the center of ${\hg}^{ext}$.

 \begin{theorem} \label{struk-alg-novs}
The vertex algebra $W(0)$ is generated by $e,f,h, j$ and
$$W(0) = U({\hg} ^{ext} ) . ({\1} \otimes {\1}) \cong V \otimes M_T (0) .$$
\end{theorem}
{\em Proof.}
Let ${\U}$ be the vertex subalgebra of ${\cV} \otimes F_{-1}$
generated by the set $\{ e, f, h, j \}$. The components of the
fields $e(z), f(z), h(z), T(z)$ span the Lie algebra ${\hg}
^{ext}$. Since the field $T(z)$ commutes with the action of
${\hg}$ we have that ${\U} = U({\hg} ^{ext}). ({\1} \otimes {\1})
\cong V \otimes M_T (0)$.

Since the operator  $H(0)$ acts trivially on the vertex algebra
${\U}$, we conclude that ${\U}  \subset W(0)$.
We shall now prove that  the vectors $e ,f,h,j$ generate $W(0)$.
First we notice that
\bea && W(0) = \oplus _{m \in {\Z} } {\cV} ^{m} \otimes F_{-1}
^{m} .\label{dek-w0} \eea

Since
$$S(-2){\1} \otimes {\1} = \tfrac{1}{4} ( e(-1) f(-1) + f(-1) e(-1)
+ {\hf} h(-1) ^{2} ). ( {\1} \otimes {\1}) - {\hf} T(-1) ^{2} {\1}
\otimes {\1},$$
we have that
$
S(-2){\1} \otimes {\1} \in {\U} $.
Therefore,
\bea {\cV} ^{com} \subset {\U}. \label{com-podalgebra}\eea

We note that the following relations hold:
 \bea
&& G^{+} (-\thf) {\1} \otimes {\1} = e(-2). ({\1} \otimes
e^{-\b}), \label{rel+} \\
&& [ e(m), ({\1} \otimes e^{-\b} )_n] = 0, \quad \mbox{for every}
\
\ m,  n \in {\Z},  \label{komut-en} \\
&& G^{-} (-\thf) {\1} \otimes {\1} = f(-2). ({\1} \otimes
e^{\b}), \label{rel-} \\
&& [ f(m), ({\1} \otimes e^{\b}) _n] = 0, \quad \mbox{for every} \
\ m,  n \in {\Z}. \label{komut-fn} \\
 && {\1} \otimes {\b}(-1) = -{\hf} ( h(-1) -2 T(-1) ) . ({\1}
\otimes {\1}) \label{rel-h} ,
 \eea
Relations (\ref{com-podalgebra}) and  (\ref{rel-h})  imply  that
\bea
&& {\cV}^{com} \otimes F_{-1} ^{0} = {\cV} ^{com} \otimes
M_{\b}(1) \subset {\U} \label{prva-pod} . \eea

We shall now prove that   ${\cV} ^{m} \otimes F_{-1} ^{m} \subset
{\U}$, where $m \in {\Z}$.
By using  the description of ${\cV}^{m}$ from
(\ref{opis-komponente}), we conclude that it is enough to show
that  elements of the form
\bea v &=&  G^{+ }({-n_1-{\thf}})  \cdots G^{+}({-n_r-{\thf}})
G^{-}({-k_1-{\thf}})  \cdots G^{-}({-k_s-{\thf}}) w \otimes
w^{(1)} \non \\ && w \in {\cV} ^{com},  w^{(1)} \in F_{-1} ^{m}, \
n_i, k_{\ell} \in {\Zp}, r-s=m, \non \eea
belong to ${\U}$.

By using (\ref{rel+}) and (\ref{komut-en}) we get that
$$ v = \sum_{i \in I} \tilde{g}_i ^{(2)} . \left( G^{-}({-k_1-{\thf}})  \cdots G^{-}({-k_s-{\thf}}) w \otimes
w_i ^{(2)} \right)  $$
for certain $w_i ^{(2)} \in F_{-1} ^{-s}$,
$ \tilde{g}_i ^{(2)} \in U ( {\C} e \otimes {\C}[t,t^{-1}] )
\subset U(\hg)$ and a finite index set $I$.
Similarly, applying (\ref{rel-}) and (\ref{komut-fn}) we get
$$ G^{-}({-k_1-{\thf}})  \cdots G^{-}({-k_s-{\thf}}) w \otimes
w_i ^{(2)} =   \sum_{\ell \in I_i} \tilde{g} _{i,\ell} ^{(3) } .(
w \otimes w_{i,\ell} ^{(3)} )$$
for certain $w_{i,\ell} ^{(3)} \in F_{-1} ^{0} = M_{\b}(1)$,
$ \tilde{g}_{i,\ell}   ^{(3)} \in U ( {\C} f \otimes
{\C}[t,t^{-1}] ) \subset U(\hg)$, and a finite index set $I_i$.
Since   $w \otimes w_{i,\ell} ^{(3)} \in {\cV} ^{com} \otimes
M_{\b}(1) $, we obtain
$$  v \in  U(\hg) .  ( {\cV} ^{com}
\otimes M_{\b}(1) ) \subset U({\hg} ^{ext}). ({\1} \otimes {\1}) =
{\U}.
$$

Therefore,
$$ {\cV}^{m} \otimes F_{-1}^{m} \subset {\U}, \quad \mbox{for
every} \ m \in {\Z}.$$
Now (\ref{dek-w0}) implies that $W(0) = {\U}$.
  \qed

\vskip 5mm

 The next Lemma follows form   Corollary 4.2 of
\cite{DM-galois}. This result will be our important tool in the
irreducibility analysis.

\begin{lemma} \label{pom-ired} Assume that $M$ is an irreducible ${\cV} \otimes F_{-1}$--module.
Then for each $0 \ne w \in M$, $M$ is spanned as a ${\cV} \otimes
F_{-1}$--module by $u_n w$, for $u \in {\cV} \otimes F_{-1}$ and
$n \in {\Z}$.
\end{lemma}

\begin{lemma} \label{pomoc-ireducibilnost}
Assume that $R$ is an irreducible $W(0)$--module. Then $R$ is an
irreducible $V$--module. In particular, $R$ is an irreducible
${\hg}$--module at the critical level.
\end{lemma}
{\em Proof.} Since $W(0) \cong V \otimes M_T(0)$, we have that $R
\cong S \otimes N$, where $S$ is an irreducible $V$--module and
$N$ is an irreducible $M_T(0)$--module. Since $M_T(0)$ is a
commutative vertex algebra, we have that $N$ is one dimensional
and that  every element $T(n)$ acts
 on $R$ as a scalar multiplication.
Therefore  $R$ is irreducible as a $V$--module. \qed

\begin{theorem} \label{osnovni-ir}
Assume that $U$ is a ${\cV}$--module such that $U$ admits the
following ${\Z}$--graduation
\bea && U = \bigoplus_{j \in {\Z} } U^{j}, \quad {\cV} ^{i} .
U^{j} \subset U ^{i+j} . \label{uvjet-grad} \eea
Then
$$ U \otimes F_{-1} = \bigoplus_{s \in {\Z} } {\mathcal L}_{s}(U), \quad  \mbox{where} \ \    {\mathcal L}_{s} (U)
 := \bigoplus_{i \in {\Z} } U^{i} \otimes F_{-1} ^{- s+i},$$
  is an  $W(0)$--module.

If $U$ is irreducible, then for every $s \in {\Z}$  $ {\mathcal
L}_{s}(U)$ is an irreducible $V$--module.
\end{theorem}
{\em Proof.} Since ${\cV} \otimes F_{-1} = \oplus_{\ell \in {\Z} }
W(\ell)$, relation (\ref{uvjet-grad}) implies that
\bea W(\ell) .  {\mathcal L}_{s}(U) \subset  {\mathcal L}_{s+
\ell}(U)\  \mbox{for every} \ \ell, s \in {\Z}. \label{pom-ir2}
\eea
This proves that  ${\mathcal L}_{s}(U)$ is a $W(0)$--module for
every $s\in {\Z}$.

Assume now that $U$ is irreducible. Then $U \otimes F_{-1}$ is an
irreducible ${\cV} \otimes F_{-1}$--module.
 Let $0 \ne v \in {\mathcal L}_{s}(U)$. Since $U\otimes
F_{-1} $ is a simple  ${\cV} \otimes F_{-1}$--module, by Lemma
\ref{pom-ired}  we get that
\bea \label{generiranje-4} &&   U \otimes F_{-1}  =
\mbox{span}_{\C} \{ u_n v \ \vert \ u \in {\cV} \otimes F_{-1}, \
n \in {\Z} \} . \eea
By using    (\ref{pom-ir2}) and (\ref{generiranje-4})  we conclude
that
\bea \label{generiranje-5} &&  {\mathcal L}_{s}(U) =
\mbox{span}_{\C} \{ u_n v \ \vert \ u \in W(0), \ n \in {\Z}  \} .
\eea
 So ${\mathcal L}_{s}(U)$
is an irreducible  $W(0)$--module. Now Lemma
\ref{pomoc-ireducibilnost} gives that ${\mathcal L}_{s}(U)$ is an
irreducible $V$--module, and therefore an irreducible
${\hg}$--module of critical level.\qed

\begin{corollary} \label{konstr-ir-mod}
Assume that $U \subset F(\chi^{+}, \chi^{-})$ is an irreducible
${\cV}$--module. Then $U \otimes F_{-1}$ is a completely reducible
$V$--module
$$ U \otimes F_{-1} = \bigoplus_{s \in {\Z} } {\mathcal L}_{s}(U)$$ and   ${\mathcal L}_{s} (U) = \{ v \in
U \otimes F_{-1} \ \vert \ H(0) v = s v \}$
 is an irreducible $V$--module. Moreover, ${\mathcal L}_{s} (U)$ is an irreducible ${\hg}$--module
  at the critical level.
\end{corollary}
{\em Proof.} The operator $J^{f}(0)$ acts semisimply on $U \subset
F(\chi^{+}, \chi^{-})$ and defines on $U$ the following graduation
$$U = \bigoplus_{j \in {\Z} } U^{j}, \quad \mbox{where} \  U^{j} =\{ u \in U \ \vert \ J^{f}(0) u = j u
\}.$$
Now Theorem \ref{osnovni-ir} implies that ${\mathcal L}_{s} (U) $
is an irreducible ${\hg}$--module at the critical level. \qed

\section{Weyl vertex algebra and irreducibility of the Wakimoto modules}
\label{wakimoto-mod}

In this section we will see that our ${\hg}$--modules include the
Wakimoto ${\hg}$--modules at the critical level defined   by using
vertex algebra $W$ associated to the Weyl algebra. As an
application, we present  a proof of irreducibility  for a family
of Wakimoto modules.

First we shall consider the simple vertex superalgebra
$\widetilde{F} \otimes F_{-1}$. The operator $H(0)$ acts semismply
on $\widetilde{F} \otimes F_{-1}$, and we have that
$$\mbox{Ker}_{\widetilde{F} \otimes F_{-1}} H(0) ={\mathcal L}_0
(\widetilde{F})$$
is a simple vertex algebra. We shall now identify this vertex
algebra.

Define: \bea \label{for-weyl} a:= \Psi^{+} (-\thf) {\1} \otimes
e^{\b}, \quad a^{*} := - \Psi^{-}(-\hf){\1} \otimes e^{-\b}, \eea
and
$$ a(z) =Y(a,z) = \sum_{n \in {\Z} } a(n) z^{-n-1}, \ \ a^{*}(z) =Y(a^{*},z) = \sum_{n \in {\Z} } a^{*}(n)
z^{-n}. $$
Then $$[a(n), a(m)] = [a^{*}(n), a^{*}(m)] = 0, \quad [a(n),
a^{*}(m)] = \delta_{n+m,0} .$$
Therefore, the components of the fields $a(z)$ and $a^{*}(z)$ span
the infinite dimensional Weyl algebra. Let $W$ be the vertex
subalgebra of $\widetilde{F}\otimes F_{-1}$ generated by $a(z)$
and $a^{*}(z)$ (cf. \cite{FMS}, \cite{efren}, \cite{a-comun-a}).

By using a similar proof to those of Theorem \ref{struk-alg-novs}
we get.

\begin{proposition} The vertex algebra ${\mathcal L}_0
(\widetilde{F})$ is generated by $a$ and $a^{*}$. Thus we have:
\bea \label{weyl-identifcation} W \cong {\mathcal L}_0
(\widetilde{F}). \eea
\end{proposition}

\vskip 5mm

 For every ${\chi} \in {\C}((z))$, the vector space $$
  \widetilde{F}_{\chi} \otimes F_{-1} \cong \widetilde{F}
\otimes F_{-1}$$ carries the structure of a ${\hg}$--module at the
critical level, and  ${\mathcal L}_0  (\widetilde{F}_{\chi} )
\subset \widetilde{F} \otimes F_{-1}$ is a ${\hg}$--submodule.
By construction, we have that as a vector space
${\mathcal L}_0  (\widetilde{F}_{\chi} )$ is isomorphic to
 ${\mathcal L}_0 (\widetilde{F} ) \cong W $. Therefore on the
 vertex algebra $W$ for every $\chi \in {\C}((z))$ exists a
 ${\hg}$--structure.

 By using the
definition of ${\hg}$--structure on $W$ and (\ref{for-weyl}) one
gets:

\bea
e(z)&=& a(z), \nonumber \\
h(z) &=& - 2 :  a^{*}(z) a(z): - {\chi}(z)           \nonumber \\
f(z) & =&   - : a^{*} (z)  ^{2} a(z) : -2 \partial_{z} a^{*} (z) -
a^{*}(z) {\chi}(z) . \nonumber
\eea

Therefore the   ${\hg}$--module  structure on  $W$ coincides with
the Wakimoto module $W_{-\chi}$ (see \cite{efren} and reference
therein).

Combining  Proposition \ref{vazna-2}   and Corollary
\ref{konstr-ir-mod} we obtain the following irreducibility
result.
\begin{theorem} \label{ireducibilnost-wakimotovih}
\item[(1)] For every ${\chi} \in {\C}((z))$, the ${\hg}$--module
${\mathcal L}_0 ({\widetilde F}_{\chi})$ is isomorphic to the
Wakimoto module $W_{-\chi}$.
\item[(2)] Assume that $${\chi}(z) = \sum_{n=-p}  ^{\infty}
{\chi}_{-n} z ^{n-1} \in {\C}((z)) \quad (p \in {\Zp})$$ satisfies
the conditions (\ref{uvjet-prvi}) and (\ref{uvjet-drugi}) of
Proposition \ref{vazna-2}. Then $W_{-\chi} \cong {\mathcal L}_0
({\widetilde F}_{\chi})$ is an irreducible ${\hg}$--module at the
critical level.
\end{theorem}

\section{Construction of irreducible highest weight  modules }
\label{konstr-najv}

 In this section we apply the results from
previous sections and obtain a  construction of all irreducible
highest weight ${\hg}$--modules $L(\l)$ of level $-2$.  It turns
out that modules $L(\l)$ are realized inside the Weyl vertex
algebra $W$, and therefore they are submodules of certain Wakimoto
modules.
% In
%generic case we obtain a new proof of irreducibility.
 By using the methods developed in \cite{A-2005}, we shall
 identify modules obtained from irreducible highest weight
modules by applying the automorphism $\pi_s$. We  will also show
the vertex superalgebra $\overline{F} \otimes F_{-1}$ is a
completely reducible module for the simple vertex algebra
$L(-2\Lambda_0)$.

\vskip 5mm

First we shall study non-generic highest weight representations.

\begin{theorem} \label{najvece-tezine-1}
\item[(i)]For every $n \in {\Zp}$ the vector space
 $\overline{F}(0,-\tfrac{n}{z})  \otimes
F_{-1} $ carries an ${\hg}$--structure uniquely determined by
\bea
e(m) &=& - \sum_{i \in {\Z} } i \Psi ^{+} (i-{\hf} ) \otimes
e^{\b} _{m -i -1} \nonumber \\
f(m) & = & - \sum_{i \in {\Z} } (i+ n) \Psi ^{-} (i-{\hf} )
\otimes
e^{-\b} _{m -i -1} \nonumber \\
h(m) & = & -2 \b (m) + n \delta_{m,0} \nonumber \\
c &=& -2 \nonumber ,
\eea where $m \in {\Z}$.
 Moreover, $\overline{F}(0,-\tfrac{n}{z})   \otimes F_{-1}$ is
a completely reducible ${{\hg}}$--module and
$$\overline{F} (0, -\tfrac{n}{z})  \otimes F_{-1}\cong \bigoplus_{s \in {\Z} } \pi_s ( L( -(2+n) \Lambda_0 + n \Lambda_1)).$$

\item[(ii)]
$$ U({\hg}). ({\1} \otimes {\1}) \cong ( L( -(2+n) \Lambda_0 + n
\Lambda_1), \ \  U({\hg}). ({\1} \otimes e^{-\b}) \cong ( L( n
\Lambda_0 - ( n + 2) \Lambda_1).$$
\end{theorem}
{\em Proof.} By using Proposition \ref{ireducibilnost-druga} we
have that for every $n \in {\Zp}$, $\overline{F}(0,-\tfrac{n}{z})$
is an irreducible ${\cV}$--module. Then Corollary
\ref{konstr-ir-mod} gives that $\overline{F}(0,-\tfrac{n}{z})
\otimes F_{-1}$ is a completely reducible ${\hg}$--module
isomorphic to $\oplus_{s \in {\Z} } {\mathcal
L}_{s}(\overline{F}(0,-\tfrac{n}{z})  )$, where ${\mathcal
L}_{s}(\overline{F}(0,-\tfrac{n}{z})  )$ is an irreducible
${\hg}$--module.
For every $s \in {\Z}$, we set
$$v_s = {\1} \otimes e^{-s\b} \in {\mathcal
L}_{s}(\overline{F}(0,-\tfrac{n}{z})). $$

Now using Lemma \ref{sing-za} one obtains that
$$ {\mathcal
L}_{s}(\overline{F}(0,-\tfrac{n}{z})  ) = U({\hg}). v_{s} =
\pi_{-s} (L(-(2+n) \Lambda_0 + n \Lambda_1)). $$
This proves (i). The second assertion follows form (i) and from
the fact that $\pi_{-1} ( L(-(2+n) \Lambda_0 + n \Lambda_1) )
\cong L(n \Lambda_0-(2+n) \Lambda_1 )$.  \qed

\vskip 5mm

When $n \in  {\N}$, then   $L( -(2+n) \Lambda_0 + n \Lambda_1)$ is
not a module for the simple vertex algebra $L(-2\Lambda_0)$. So
Theorem  \ref{najvece-tezine-1} can be applied only in the
framework of $N(-2\Lambda_0)$--modules. But  when $n=0$ we have
that ${\cV}$--module $\overline{F}=\overline{F}(0,0)$ is a simple
vertex superalgebra (see Proposition \ref{simple-quotient}), and
we have the following realization of $L(-2\Lambda_0)$--modules.
\begin{corollary}
The simple vertex algebra $L(-2\Lambda_0)$ is a subalgebra of the
vertex superalgebra $ \overline{F} \otimes F_{-1}$, and we have
the following decomposition of $L(-2\Lambda_0)$--modules:
$$\overline{F} \otimes F_{-1} \cong \oplus_{s \in {\Z} } \pi_s
(L(-2\Lambda_0)) .$$
\end{corollary}
{\em Proof.} Since $\overline{F} \otimes F_{-1}$ is a vertex
superalgebra we have that $$\mbox{Ker}_{\overline{F} \otimes
F_{-1} } H(0) \cong {\mathcal L}_0 (\overline{F} ) \cong L(-2
\Lambda_0) $$ is a vertex subalgebra of $\overline{F} \otimes
F_{-1}$. Now the statement  follows from Theorem
\ref{najvece-tezine-1}. \qed \vskip 5mm

 By using  Corollary \ref{ireducibilnost-treca} and the
proof similar to those of Theorem \ref{najvece-tezine-1} one
obtains the following result.

\begin{theorem} \label{najvece-tezine-2}
For every ${\l} \in \{-1 \} \cup {\C} \setminus {\Z}$, the vector
space $\widetilde{F}_{-\tfrac{\l}{z} } \otimes F_{-1}$ carries an
${\hg}$--structure uniquely determined by :
\bea
e(m) &=& - \sum_{i \in {\Z} } i \Psi ^{+} (i-{\hf} ) \otimes
e^{\b} _{m -i -1} \nonumber \\
f(m) & = & - \sum_{i \in {\Z} } (i+{\l}) \Psi ^{-} (i-{\hf} )
\otimes
e^{-\b} _{m -i -1} \nonumber \\
h(m) & = & -2 \b (m) + {\l} \delta_{m,0} \nonumber \\
c &=& -2 ,\nonumber \eea
where $m  \in  {\Z}$.
 Moreover, $\widetilde{F}_{-\tfrac{\l}{z}
}\otimes F_{-1}$ is a completely reducible ${\hg}$--module and
$$\widetilde{F}_{-\tfrac{\l}{z} } \otimes F_{-1}\cong \bigoplus_{s \in {\Z} } \pi_s ( L( -(2+{\l}) \Lambda_0 +
{\l} \Lambda_1)).$$
\end{theorem}

\section{Realization of irreducible  modules on the vertex algebra $\Pi(0)$}
\label{konstr-relaxed}

So far we studied the irreducible ${\hg}$--modules realized on the
Weyl vertex algebra. In this section we shall see that there
exists a family of irreducible ${\hg}$--modules realized on a
larger vector space. In order to construct new irreducible
representations, we shall study the vertex algebra
$\Pi(0)={\mathcal L}_0(F)$ which contains the Weyl vertex algebra
$W$ as a subalgebra.

We shall also identify the  irreducible ${\hg}$--module which are
${\Zp}$--graded but don't belong to the category $\mathcal{O}$.
These modules are irreducible quotients of relaxed Verma modules
studied in \cite{FST}.

First we recall that the boson-fermion correspondence gives that
the fermionic vertex superalgebra $F$ is isomorphic to the lattice
vertex superalgebra $V_{\Z}$. Therefore,
$$ F \otimes F_{-1} \cong V_L = M(1) \otimes {\C}[L],$$

where the lattice $V_L$ is the lattice vertex superalgebra
associated  to the lattice

$$L= {\Z} {\a} + {\Z} {\b} , \quad \la \a , \a \ra = - \la \b , \b
\ra = 1, \ \la \a , \b \ra = 0.$$
(As usual, $M(1)$ is a level one irreducible module for the
Heisenberg algebra $\hat{\h}_{\Z}$ associated to the   abelian
algebra ${\h}= L \otimes_{\Z} {\C}$ and ${\C}[L]$ is the group
algebra with  generators $e^{\a}$ and $e^{\b}$.)

The operator $H$ from  (\ref{def-operatora-h}) coincides with
${\a}(-1) + {\b}(-1)$.

 We conclude that the vertex algebra $\Pi (0)=\mbox{Ker}_{F\otimes F_{-1}} ({\a} + {\b})(0) ={\mathcal L}_0 (F)$ is
  isomorphic to the vertex algebra
 \bea \label{izomorfizam-resetka}&& \Pi(0) \cong M (1) \otimes {\C} [{\Z}(\a + \b)] . \eea

\begin{remark}
By using different methods, the vertex algebra $\Pi(0)$ was  also
studied by E. Frenkel in \cite{efren} and by  S. Berman, C. Dong
and S. Tan in  \cite{BDT}.
\end{remark}

Now we shall apply the results from Sections \ref{konstr-a-ired}
and  \ref{mappings}. The following result gives a construction of
a large class of irreducible ${\hg}$--modules on the vertex
algebra $\Pi(0)$.

Combining Proposition \ref{vazna-1} and Corollary
\ref{konstr-ir-mod} we obtain the following result.

\begin{theorem} \label{ired-chi-pi}
Assume that ${\l} \in {\C}\setminus {\Z}$ and that $${\chi}(z) =
\sum_{n=-p} ^{\infty} {\chi}_{-n} z ^{n-1} \in {\C}((z)) \quad (p
\in {\Zp})
$$ satisfies the conditions (\ref{uvjet-pr}) and (\ref{uvjet-dr})
of Proposition \ref{vazna-1}. Then on the vertex algebra $\Pi(0)=
{\mathcal L}_0(F)$ exists the structure of an irreducible
${\hg}$--module isomorphic to ${\mathcal L}_0 (F(\tfrac{\l}{z},
{\chi}))$.
\end{theorem}

\vskip 5mm

Define the following Virasoro vector in $ \Pi(0) \subset F\otimes
F_{-1}$:
$$\omega= \omega^{(f)} \otimes {\1} - {\hf} {\1} \otimes {\b}(-1)
^{2}= {\hf} ( {\a} (-1) ^{2} - {\b}(-1) ^{2}) . ({\1} \otimes
{\1}) .$$
Let $L(z)=Y(\omega,z) = \sum_{n \in {\Z} } L(n) z^{-n-2}$. Then
the components of the field $L(z)$ satisfies the commutation
relations for the Virasoro algebra with central charge $c=2$.
Moreover,  $L(0)$ acts semisimply on $F\otimes F_{-1}$ with
half-integer eigenvalues, and it defines a ${\Zp}$--graduation on
the vertex algebra $\Pi(0)$ :
\bea \label{grad-pi}  \Pi(0) = \bigoplus_{m \in {\Zp} } \Pi(0)_m ,
\quad \Pi(0)_m = \{ v \in \Pi(0) \ \vert \ L(0) v = m v \}.\eea
Let
$\mbox{ch}_{\ \Pi(0)} (q,z) = \mbox{tr} \  q^{L(0)} z^{- 2 \b(0)}.
$

By using relation (\ref{izomorfizam-resetka}) and the  properties
of the $\delta$-function  one can easily show the following
result.
\begin{proposition} We have:
\bea \label{ch-pi} \mbox{ch}_{\ \Pi(0)}(q,z) &=&  \delta (z ^{2})
\prod_{n =1} ^{\infty} (1-q^{n}) ^{-2}  \nonumber \\
&=& \delta (z ^{2}) \prod_{n =1} ^{\infty} (1-q^{n}z ^{2}) ^{-1}
(1-q^{n} z^{-2}) ^{-1}  .\eea \end{proposition}

We are now interested in ${\hg}$--modules from Theorem
\ref{ired-chi-pi} such that the ${\hg}$--action is compatible with
the $L(0)$--graduation on the vertex algebra $\Pi(0)$.   But the
action is compatible with the graduation if and only if ${\chi}(z)
= \frac{\mu}{z} $ for certain ${\mu} \in {\C}$. Therefore we
should consider the irreducible ${\cV}$--module
$F(\tfrac{\l}{z},\tfrac{\mu}{z})$, where ${\l}, {\mu} \in {\C}
\setminus {\Z}$.  The module $F(\tfrac{\l}{z},\tfrac{\mu}{z}) $
has the  simple structure as an ${\cA}$--module. In fact, the
action of the Lie superalgebra ${\cA}$ is (up to scalar factor)
the same as the action of the Clifford algebra $CL$ on $F$. When
we apply Corollary \ref{konstr-ir-mod}, we get a family of
irreducible ${\hg}$--modules ${\mathcal L}_s
(F(\tfrac{\l}{z},\tfrac{\mu}{z}))$ at the critical level.

Now we want to identify these  irreducible ${\hg}$--modules.

 For every  $s \in {\Z}$, we  define a family of vectors
  $w_j ^{(s)}  \in {\mathcal L}_s
  (F(\tfrac{\l}{z},\tfrac{\mu}{z}))$, by
\bea
w^{(s)}_0 := && {\1} \otimes e^{-s{\b}}, \nonumber \\
w^{(s)}_{j}:=&&   \Psi^{+}(-j+{\hf}) \cdots
\Psi^{+}(-\tfrac{1}{2}){\1}
\otimes e^{(j - s)\b}
%\in  {\mathcal L}_s (F(\tfrac{\l}{z},\tfrac{\mu}{z}))
\quad (j \in {\N}), \nonumber \\
w^{(s)}_{-j}:=  && \Psi^{-}(-j+{\hf}) \cdots
\Psi^{-}(-\tfrac{1}{2}){\1} \otimes e^{-(j+s)\b} \quad (j \in
{\N}).
% \in  {\mathcal L}_s (F(\tfrac{\l}{z},\tfrac{\mu}{z}))
\nonumber \eea
 By using a direct calculation, one can prove the
following lemma.

\begin{lemma} \label{def-top-level}
Assume that   $s,j \in {\Z}$ and $n \in {\Zp}$. Then we have
%
%$$ e(n+s). w ^{(s)} _{j} =   f(n - s). w ^{(s)} _{j} =  h (n). w ^{(s)} _{j} =  0  \mbox{for every} \  n \in {\N}.$$
%
\bea  && e(n - s) . w^{(s)}_{j} = \delta_{n,0}  (\l +j)
w^{(s)}_{j+1}, \nonumber
\\
&& h(n). w^{(s)}_j =  \delta_{n,0} ( 2 j - 2 s + \l - \mu)
w^{(s)}_j, \nonumber
\\ &&  f(n+s). w^{(s)}_{j} =  \delta_{n,0} (\mu - j) w^{(s)}_{j-1}. \nonumber \eea
\end{lemma}

\vskip 5mm

Let us first consider the case $s=0$. Define:
$$ E_{\l,\mu} := {\mathcal L}_0 (F(\tfrac{\l}{z},\tfrac{\mu}{z})).$$ By construction, $E_{\l,\mu} \cong \Pi(0)$ as a vector space.
  Introduce the graduation operator $L(0)$ on
the ${\hg}$--module $E_{\l,\mu}$ by using the vertex algebra
graduation (\ref{grad-pi}) on $\Pi(0)$. So let
$$E_{\l,\mu} =\bigoplus_{m\in {\Zp} } E_{\l,\mu} (m), \ \ E_{\l,\mu}(m) = \{ v \in E_{\l, \mu} \ \vert \ L(0) v = m v \}.
$$
  Since $[L(0), x(n)] =-n x(n)$ for every $x\in {\g}$,  we have
  that ${\hg}$--action on $E_{\l,\mu}$ is compatible with the
  graduation. In other words, $E_{\l,\mu}$ is an ${\hg} \oplus
  {\C}L(0)$--module.
Lemma \ref{def-top-level} shows that  the top level
$$E_{\l,\mu}(0) = \mbox{span}_{\C} \{ w^{(0)}_j \ \vert \ j \in {\Z} \}$$
is an irreducible $U({\g})$--module which is neither highest nor
lowest weight with respect to ${\g}$.

Next we consider general case. By using Lemma \ref{def-top-level}
we have that $\pi_s \left({\mathcal L}_s
(F(\tfrac{\l}{z},\tfrac{\mu}{z}))\right) $ is an irreducible
${\Zp}$--graded ${\hg}$--module whose top level is isomorphic to
$E_{\l,\mu}(0)$. This proves that $\pi_s \left({\mathcal L}_s
(F(\tfrac{\l}{z},\tfrac{\mu}{z}))\right) \cong E_{\l,\mu}$.
Therefore, $${\mathcal L}_s (F(\tfrac{\l}{z},\tfrac{\mu}{z}))\cong
\pi_{-s} (E_{\l,\mu}).$$

In this way we have proved the following result.

 \begin{theorem} \label{najvece-tezine-3}
For every $\l, \mu \in {\C} \setminus {\Z}$, the vector space
$F(\tfrac{\l}{z},\tfrac{\mu}{z}) \otimes F_{-1} \cong F \otimes
F_{-1}$ carries an ${\hg}$--structure uniquely determined by :
\bea
e(m) &=&  \sum_{i \in {\Z} } (\l -i) \Psi ^{+} (i-{\hf} ) \otimes
e^{\b} _{m -i -1} \nonumber \\
f(m) & = &  \sum_{i \in {\Z} } (\mu -i) \Psi ^{-} (i-{\hf} )
\otimes
e^{-\b} _{m -i -1} \nonumber \\
h(m) & = & -2 \b (m) +   (\l-\mu ) \delta_{m,0} \nonumber \\
c &=& -2 ,\nonumber \eea
where $m \in {\Z}$.
Moreover, $F(\tfrac{\l}{z},\tfrac{\mu}{z}) \otimes F_{-1}$ is a
completely reducible ${\hg}$--module and
$$F(\tfrac{\l}{z},\tfrac{\mu}{z}) \otimes F_{-1}\cong \bigoplus_{s \in {\Z} }  \pi_{s}
(E_{\l,\mu}).$$
\end{theorem}

Although modules  $E_{\l,\mu}$ don't belong to the category
$\mathcal{O}$, one can investigate their characters.  The
operators $h(0)$ and $L(0)$ acts semisimply on $E_{\l,\mu}$ with
finite-dimensional common eigenspaces. Since we need the degree
operator $L(0)$ we shall consider $E_{\l,\mu}$ as a module for the
Lie algebra ${\hg} \oplus {\C}L(0)$.

\begin{corollary} \label{relaxed-char}
For every $\l,\mu \in {\C} \setminus {\Z}$, the vertex algebra
$\Pi(0)$ carries the structure of an irreducible ${\Zp}$--graded
${\hg} \oplus {\C}L(0)$--module isomorphic to $E_{\l,\mu}$. We
have the following character formulae:
\bea \mbox{ch}_{E_{\l,\mu}} (q,z) &=& \mbox{tr} \ q^{L (0) } z^{
h(0)} =z^{ \l-\mu}   \delta (z ^{2}) \prod_{n =1} ^{\infty}
(1-q^{n}z ^{2})
^{-1} (1-q^{n} z^{-2}) ^{-1} \nonumber \\
&=& z^{ \l-\mu} \delta (z ^{2}) \prod_{n =1} ^{\infty} (1-q^{n})
^{-2}. \nonumber \eea
\end{corollary}
{\em Proof.} We already sow that $E_{\l,\mu}$ is an irreducible
${\Zp}$--graded ${\hg} \oplus {\C} L(0)$--module.   By using
(\ref{ch-pi}) we get
$$\mbox{ch}_{E_{\l,\mu} } (q,z) = \mbox{tr} \ q^{L (0)}
z^{ h (0)} = z^{ \l-\mu} \mbox{ch}_{\Pi(0)}(q,z) =z^{{\l}-{\mu}}
\delta (z ^{2}) \prod_{n =1} ^{\infty} (1-q^{n}) ^{-2}.$$
This proves the character  formula. \qed

\begin{remark}
 Modules $E_{\l,\mu}$ are irreducible quotients of certain relaxed Verma
 modules which are introduced and studied  in \cite{FST}. In our
 terminology, the relaxed Verma modules are $N(-2, E_{\l,\mu} (0))$ and
 they have
 the following character
$$ \mbox{ch}_{N(-2, E_{{\l},\mu} (0))}(q,z) = z^{ \l-\mu}\delta (z^{2}) \prod_{n =1} ^{\infty} (1-q^{n})
^{-3}.$$
In Corollary \ref{relaxed-char} we calculate the characters of
irreducible quotients of relaxed Verma modules by  using our
explicit realization. So our methods  don't use  structure theory
of relaxed Verma modules.
\end{remark}


\begin{thebibliography}{Lia}



\bibitem [A1] {A3}  D. Adamovi\' c, Representations of the $N=2$
superconformal vertex algebra, Internat. Math. Res. Notices  {\bf
2} (1999) 61-79.



\bibitem [A2] {a-comun-a} D. Adamovi\' c, Representations of the vertex
algebra ${\mathcal W}_{ 1 + \infty}$ with a negative integer
central charge, Comm. Algebra 29 (2001) no. 7, 3153-3166.

\bibitem [A3]{A-2005} D. Adamovi\' c, A construction of admissible
$A_1^{(1)}$--modules of level $-\tfrac{4}{3}$, J. Pure  Appl.
Algebra {\bf 196} (2005) 119-134.

\bibitem [AM]{AM} D. Adamovi\'{c} and A. Milas,
Vertex operator algebras associated to the modular invariant
representations for $A_1^{(1)}$, Math. Res. Lett. {\bf 2} (1995)
563-575.

\bibitem [BDT]{BDT} S. Berman, C. Dong and S. Tan, Representations
of a class of lattice type vertex algebras,  J. Pure  Appl.
Algebra {\bf 176} (2002) 27-47.



\bibitem [DL]{DL} C.  Dong and J. Lepowsky,
 Generalized vertex algebras and relative vertex operators,
  Birkh\"auser,  Boston, 1993.



 \bibitem[DM]{DM-galois} C. Dong and G. Mason, On quantum Galois
 theory,  Duke Math. J. {\bf 86} (1997) 305-321.



\bibitem [EG]{EG} W. Eholzer and M. R. Gaberdiel, Unitarity of rational
$N=2$ superconformal theories, Comm. Math. Phys. {\bf 186} (1997)
61-85.

\bibitem [FB]{FB} E. Frenkel and D. Ben-Zvi, Vertex algebras and
algebraic curves, Mathematical Surveys and Monographs; no. 88,
AMS, 2001.

\bibitem[F]{efren} E. Frenkel, Lectures on Wakimoto modules,
opers and the center at the critical level, Adv. Math 195 (2005)
297-404.

\bibitem[FF]{FF1} E. Frenkel and  B. Feigin, Affine Kac-Moody algebras and semi-infinite flag manifolds,
 Comm. Math. Phys. {\bf 128} (1990) 161-189.


\bibitem [FHL]{FHL}
I. B. Frenkel, Y.-Z. Huang and J. Lepowsky, On axiomatic
approaches to vertex operator algebras and modules, Mem. Amer.
Math. Soc. {\bf 104}, 1993.

\bibitem [FLM]{FLM}
I. B. Frenkel, J. Lepowsky and A. Meurman,   Vertex Operator
Algebras and the Monster,   Pure and Applied Math., Vol. {\bf
134}, Academic Press, New York, 1988.

\bibitem[FMS]{FMS} D. Friedan, E. Martinec and S. Shenker,
Conformal invariance, supersymmetry and string theory, Nuclear
Phys. B {\bf 271} (1986) 93-165.



\bibitem [FST]{FST} B. L. Feigin, A. M. Semikhatov and  I. Yu. Tipunin,
Equivalence between chain categories of representations of affine
$sl(2)$ and $N=2$ superconformal algebras, J. Math.  Phys. {\bf
39} (1998), 3865-390.

\bibitem[FZ]{FZ}
I. B. Frenkel and Y.  Zhu, Vertex operator algebras associated to
representations of affine and Virasoro algebras,  Duke Math. J.
{\bf 66} (1992),  123-168.


\bibitem [K1]{K-b} V. Kac, Infinite dimensional Lie algebras, Third edition, Cambridge Univ. Press, Cambridge, 1990.


\bibitem [K2]{K}  V. Kac,   Vertex Algebras for Beginners, University
Lecture Series, Second Edition,   Amer. Math. Soc., 1998, Vol. 10.

\bibitem [KK]{KK} V. Kac and D. Kazhdan, Structure of
representations with highest weight of infinite dimensional Lie
algebras, Adv. Math. {\bf 34} (1979) 97-108.


\bibitem [KS]{KS} Y. Kazama and H. Suzuki, New N=2 superconformal
field theories and superstring compactifications, Nucl. Phys. B
{\bf 321} (1989), 232-268.






\bibitem [Li1]{Li}
H. Li,  Local systems of vertex operators, vertex superalgebras
and modules,  J. Pure  Appl. Algebra {\bf 109} (1996), 143-195.



\bibitem[Li2]{Li5} H. Li, The phyisical superselection principle
in vertex operator algebra theory, J. Algebra {\bf 196} (1997)
436-457.

\bibitem [LL]{LL} J. Lepowsky and H. Li, Introduction to vertex
operator algebras and their representations, Progress in Math.,
Vol. 227,  Birkh\"auser, Boston, 2004.


\bibitem[MP]{MP}
A. Meurman and M. Primc, Annihilating fields of standard modules
of $\tilde{sl}(2,{\C})$ and combinatorial identities, Mem. Amer.
Math. Soc. {\bf 652} 1999.


\bibitem[S]{Scz} M. Szczesny, Wakimoto modules for twisted affine
Lie algebras,  Math. Res. Lett. {\bf 9} (2002), no.4, 433-448.

\bibitem[W]{W-mod} M. Wakimoto, Fock representations of affine Lie
algebra $A_1^{(1)}$, Comm. Math. Phys. 104 (1986) 605-609.

%\bibitem [Wn]{W} W.Wang, Rationality of Virasoro Vertex operator
%algebras, Internat. Math. Res. Notices, Vol {\bf 71}, No.1
%(1993),197-211.




\bibitem[Z]{Z}
Y. Zhu, Vertex operator algebras, elliptic functions and modular
forms, Ph. D. thesis, Yale University, 1990; Modular invariance of
characters of vertex operator algebras,   J. Amer. Math. Soc. {\bf
9} (1996), 237-302.

\end{thebibliography}
\end{document}